\documentclass[10pt,a4paper]{amsart}
\usepackage{amsfonts}
\usepackage{amsthm}
\usepackage{amssymb}
\usepackage{amsmath}
\usepackage{amscd}

\theoremstyle{plain}
\newtheorem{thm}{Theorem}[section]
\newtheorem{lem}[thm]{Lemma}

\newtheorem{prop}[thm]{Proposition}
\theoremstyle{definition}
\newtheorem{defn}[thm]{Definition}

\newtheorem{ques}[thm]{Question}
\newtheorem{prob}[thm]{Problem}
\newtheorem{claim}[thm]{Claim}

\theoremstyle{remark}

\newcommand{\qq}{\mathbb{Q}}
\newcommand{\zz}{\mathbb{Z}}
\newcommand{\rr}{\mathbb{R}}
\renewcommand{\b}{\mathfrak{b}} 
\renewcommand{\d}{\mathfrak{d}} 

\newcommand{\hia}{\mathfrak{hia}}
\DeclareMathOperator{\ima}{Im}

\newenvironment{proofclaim}[1][Proof of the claim]{\textbf{#1.} }{\hfill \rule{0.5em}{0.5em}}

\title[Hereditary Interval Algebras]{Hereditary Interval Algebras and Cardinal Characteristics of the Continuum}
\author[Michael Hru\v{s}\'ak et al.]{Michael Hru\v{s}\'ak}
\address{Centro de Ciencias Matem\'aticas\\ Universidad Nacional Aut\'onoma de M\'exico\\ Campus Morelia\\Morelia, Michoac\'an\\ M\'exico 58089}
\email{michael@matmor.unam.mx}
\urladdr{http://www.matmor.unam.mx/~michael}

\author[]{Carlos Azarel Mart\'inez-Ranero}
\address{Departamento de Matem\'atica\\Facultad de Ciencias F\'isiscas y Matem\'aticas, Universidad de Concepci\'on\\Casilla 160-C, Concep\-ci\'on, Chile.}
\curraddr{}
\email{cmartinezr@udec.cl}
\author[]{Ulises Ariet Ramos-Garc\'ia}
\address{Centro de Ciencias Matem\'aticas\\ Universidad Nacional Aut\'onoma de M\'exico\\ Campus Morelia\\Morelia, Michoac\'an\\ M\'exico 58089}
\email{ariet@matmor.unam.mx}
\date{}
\thanks{The research of the first author was supported  by a PAPIIT grant IN100317 and CONACyT grant 285130.  
The  second named author was  supported by a Proyecto FONDECYT Iniciaci\'on No. 11130490 and by a Proyecto VRID-Enlace No. 218.015.022-1.0. The  third named author was partially supported by a Proyecto FONDECYT Iniciaci\'on No. 11130490, and by the PAPIIT grants  IA100517 and  IN104419.}

\begin{document}
\maketitle
\begin{abstract}
An  interval algebra is a Boolean algebra which is isomorphic to the algebra of finite unions of half-open intervals, of a linearly ordered set. An interval algebra is hereditary if every subalgebra is an interval algebra. We answer a question of M. Bekkali and S. Todor\v{c}evi\'c, by showing that  it is consistent that every $\sigma$-centered interval algebra of size $\mathfrak{b}$ is hereditary. We also show that there is, in ZFC, an hereditary interval algebra of cardinality $\aleph_1.$ 
\end{abstract}


\section{Introduction}
An $interval\ algebra $ is a Boolean algebra which is isomorphic to the algebra of finite unions of half-open intervals, of a linearly ordered set (eq., generated by a chain in the Boolean algebra ordering) and  a $subinterval \ algebra$ is a Boolean algebra which can be embedded into an interval algebra. These algebras were considered long time ago by A. Mostowski and A. Tarski \cite{MT}, who proved that every countable Boolean algebra is an interval algebra, and have been an active topic of research ever since.  For basics   results  about them,  the reader can consult  section  15 of \cite{MB} and the references therein.

\medskip

It is easy to see that there are subinterval algebras  that are not  interval algebras, the simplest example being the algebra of finite and cofinite subsets of $\omega_1$. In  this note we investigate the class of {\it hereditary interval algebras}, i.e.,   Boolean algebras with the property that all its subalgebras are interval algebras. In view of the above  mentioned  result of A. Mostowski and A. Tarski, this class contains all countable Boolean algebras.   The question of the  existence of an uncountable hereditary interval algebra, is more complicated. The finite-cofinite algebra over $\omega_1$, shows that every hereditary interval algebra must satisfy the countable chain condition.  Moreover,   M. Bekkali and S. Todor\v{c}evi\'c (\cite{BeTo}) went further and showed that every hereditary interval algebra is $\sigma$-centered. However, not all $\sigma$-centered interval algebras are hereditary. In fact, the clopen algebra of  the Alexandroff's double arrow space is not hereditary (see \cite{Nikiel}, \cite{Od}).  Taking into account the above results there is a natural cardinal invariant associated to the class of hereditary interval algebras;

\begin{defn}
Let $\mathfrak{hia}$ be  $$\min\{\,|B|\colon B\ \textrm{is\ a\ non-hereditary\ interval\ $\sigma$-centered\ algebra\,}\}.$$
\end{defn} 
This cardinal was, implicitly, investigated  by M. Bekkali and S. Todor\v{c}evi\'c in  \cite{BeTo} where they proved that:

\begin{thm}
Every $\sigma$-centered subinterval algebra $B$ of cardinality less than $\mathfrak b$ is an interval algebra.
\end{thm}
Thus, they proved that $\mathfrak b$ is a lower bound for $\mathfrak{hia}$, and moreover, it shows that the existence of uncountable hereditary interval algebras is consistent with ZFC. They also  asked the following question:
\begin{ques}\label{Q:Principal}
Are the cardinal invariants  $\mathfrak b$ and $\mathfrak{hia}$ equal?
\end{ques}

In this paper we answer Question \ref{Q:Principal} in the negative. The main ingredient of the proof is  a combinatorial reformulation of the cardinal invariant $\hia$. This reformulation comes from a careful analysis of the arguments in \cite{BeTo} and it  is captured in the following notion.  


\begin{defn} Given $A\subseteq 2^{<\omega}$, we say that  $A$ is $adequate$ if there is a function $\sigma\colon 2^{<\omega}\rightarrow 2$ such that for every $s\in 2^{<\omega}$ there are infinitely many $n$ such that $s^\frown\sigma(s)^n \in A$, where  $\sigma(s)^n$ denotes the constant function of length $n$.
 \end{defn}
The following are the main results  of the paper.

\begin{thm}\label{main1}
The cardinal $\mathfrak{hia}$ is equal to the minimal cardinality of a subset $X$ of the irrationals numbers in the Cantor set that intersect the $G_\delta$-closure of every adequate set.

\end{thm}
\begin{thm}\label{main2} Given any uncountable regular cardinal $\kappa$, it is consistent that $\hia=\kappa$ and $\mathfrak b=\aleph_1$. 

\end{thm}
It is worth mentioning that one of the main obstacles in proving this results comes from the fact that $\b\leq \hia\leq \b_2$, where $\b_2$ is a mild strengthening of $\b$ (see \ref{b2}). So roughly speaking  $\hia$ is almost $\b$.  We also complement the above results with the following.
\begin{thm}\label{main3}
There is, in ZFC, a hereditary interval algebra of cardinality $\aleph_1$.
  \end{thm}

The paper is organized as follows. In  Sections 2 and 3 will give a proof of Theorem \ref{main1}. Section 4  will give a proof of Theorem \ref{main2}. In
Section 5 will give a proof of Theorem \ref{main3}. Finally, in Section 6 we state some open problems.

\medskip

The notation and terminology in this paper is fairly standard. We
will use \cite{Jech} and \cite{Kunen} as general references for set theory,  \cite{Stevo} as a
reference for linear orders and \cite{MB} as a reference for Boolean algebras.

 \section{Upper bound for the hereditary interval number.}

In this section we begin the proof of Theorem \ref{main1}. We will  introduce a cardinal invariant $\mu$, which is a variation of the (un)bounding number $\b$, and we will show that it is equal to $\hia$. In the current section we will show that    $\hia\leq \mu$,  deferring the proof of the inequality   $\mu\leq \hia$ to the next section.  

Before proceeding any further we need to recall a few   definitions  and fix some notation. 

\begin{defn}
Let $(L,\leq_L)$ be a linearly ordered set. A subset $D$ of $L$ is called order-dense if for any $x<_L y\in L$ there is a $d\in D$ so that $x\leq_L d\leq_L y$. We say that $L$ is order-separable if there is a countable order-dense subset.
\end{defn}
\begin{defn}
Let $(L,\leq_L)$ be a linearly ordered set, $D\subseteq L $ and $f\in 2^D$. Let  $(L,\leq_L,\tau_f)$ denote the  generalised order space whose underlying set  is $L,$ and whose  topology is  generated by the subbase 

\begin{align*}
 \{\,(-\infty, x)\colon x\in L\,\} &\cup  \{\,(x,\infty)\colon  x\in L\,\}\cup \{\,(-\infty, x]\colon  x\in f^{-1}(0)\,\} \\{}& \cup\{\,[x,\infty)\colon  x\in f^{-1}(1)\,\}.
\end{align*}
\end{defn}

A \emph{Boolean} space is a compact, Hausdorff, zero dimensional topological space. Recall that  Stone duality  gives us an equivalence between the category of Boolean algebras and the category Boolean spaces, under this equivalence subalgebras correspond to quotients.  

\begin{defn}
Let $\sim$ be an equivalence relation on a Boolean space $X$. A subset $M$ of $X$ is $\sim$-saturated  
if $x\in M$ and $x\sim x'$ implies $x'\in M$. We say that $\sim$ is a Boolean equivalence relation if the subalgebra
$$B_\sim=\{\,b\in Clop(X)\colon b \textrm{ is $\sim$-saturated}\}$$ of $Clop(X)$  separates the equivalence classes.
\end{defn}
We now introduce the cardinal invariant that will be the main object of study in this section.
\begin{defn}
Let $\mu$ be equal to the minimal cardinality of an order-separable linear order $L,$ for which there exists a countable subset $ D\subseteq L$, such that $  D$ is not $G_\delta$ in $(L,\leq_L,\tau_f) $ for any function  $f\in 2^D$.
\end{defn}

\begin{prop}
The cardinal $\mu$ is equal  to  $$\nu:=\min\{\,|X|\colon \mathbb Q\subseteq X\subseteq \mathbb R, \ \mathbb Q \ \emph{is not} \ G_\delta\ \emph{in}\  (X,\leq,\tau_f)\  \emph{for any}  \ f\in 2^D\,\}.$$
\end{prop}
\proof
First notice  that $\mu\leq \nu$. Thus, it suffices to show that   $\nu\leq \mu$. To this end, we shall show that if $\kappa<\nu, $ then $\kappa<\mu$. Let $L$ be an order-separable linearly ordered set of cardinality less than $\nu$ and $D\in[L]^\omega$ be given.    Define $\hat L$ to be the linearly ordered set obtained from $L$ by inserting a copy of the rationals $\mathbb Q$ between any jump of $L$ and to  the left (right) of  the minimum  (maximum) if they exist. Let   $\hat D$  be equal to $D$ union all   added copies of $\mathbb Q$. It is easy to verify that $\hat L$ is a densely ordered set, and hence, $\hat D$ is a countable dense linear order without end-points. By Cantor's theorem there is an order preserving bijection from $\hat D$ onto $\mathbb Q$, and moreover, this can be extended to an order preserving   injection $\varphi$ from $\hat L$ into $\mathbb R$.   Set $X:=\varphi''[L]$. Since $|X|<\nu,$ there is a function    $f\in 2^\qq$ and a $G_\delta$-set $G$ in $(X,\leq, \tau_f)$  such that $G\cap X=\mathbb Q$. It follows that $\hat G=\varphi^{-1}(G)$, is a $G_\delta$-set in $(\hat L,\leq_{\hat L}, \tau_{f\circ \varphi})$ so that  $\hat G\cap \hat L=\hat D$ and $\hat G\cap L=D$. Thus, if the topology on $(L,\leq_L, \tau_g)$, where $g:= f\circ\varphi\restriction D$, coincide with the subspace topology of $(\hat L,\leq_{\hat L}, \tau_{f\circ \varphi})$, then it  follows that $\hat G\cap L$ is $G_\delta$ on $(L,\leq_L, \tau_g)$.  This fact follows from the observation that for any open ray (closed ray), say $(q,\infty)_{\hat L}$ with $q\in\hat L\setminus L$,  the set  $ (q,\infty)_{\hat L}\cap L=(a,\infty)_L$, where $q$ is any element that belongs to the jump $a<q<b$. The other cases are analogous. \endproof

In the next proposition we show that $\mu\geq \mathfrak{hia}$.  This is a generalization of Theorem 4.2 of     \cite{Nikiel}, since we need some facts from its proof, we will reprove the part relevant to us and also borrow some of its notation. 
\begin{prop}
There is a $\sigma$-centered interval algebra of cardinality $\mu$ which is not hereditary.
\end{prop}

\proof 
We will construct the desired Boolean algebra via Stone duality, i.e., we shall construct a topologically separable,  Boolean ordered  space of weight $\mu$ which admits a Boolean quotient which is not orderable.    

\medskip

Let $Q:=\mathbb Q\cap (0,1)\subseteq X\subseteq (0,1)$ be a set of cardinality $\mu$ such that   $ Q$ is not $G_\delta$ in $(X,\leq, \tau_f)$ for any  $f\in  2^Q$, and moreover, we may also assume that $X\setminus Q$ is dense in $(0,1)$.  Define $A(X)$ to be the topological space whose underlying set is $[0,1]\times\{0\}\cup X\times \{1\}\cup Q\times \{2\}$ and its topology is given by the order topology induced from the lexicographical order    (say, $<_\ell$). Note that $A(X)$  is a topologically separable, Boolean ordered space of weight $\mu$. Let $\sim$  be the equivalence relation consisting of the equality relation and all the pairs $\{(q,0),(q,2)\}$ for $q\in Q$. It is easy to verify, using   the assumption that $X\setminus Q$ is dense in $(0,1)$, that $\sim$ is a Boolean equivalence relation.  Define $\pi:A(X)\rightarrow A(X)/\sim$ to be the quotient map and give $Z:=A(X)/\sim$  the quotient topology, and let $\pi_1$ denote the projection map to the first coordinate  $\pi_1:Z\rightarrow [0,1]$.  
 
\medskip 
 
 Now, it suffices to prove that $Z$ is not orderable. Suppose for a contradiction that  $\prec$ is a linear order on  $Z$ which induces the quotient topology.  
 
 \medskip
 
 Let $Y=Z\setminus Q\times\{1\}$, and let $y_0, z_0 (y_1, z_1)$ be the $\prec$-minimum ($\prec$-maximum) elements of $Y$ and $Z$, respectively.   
 Since $Y$ is compact, it follows that for every $q\in Q$, either $z_0\preceq (q,1)\prec y_0$,   or $y_1\prec (q,1)\preceq z_1$ or $a\prec (q,1)\prec b$, where $a, b $ is a jump in $Y$, i.e.,   $(a,b)_\prec\cap Y=\emptyset$. This will be denoted by  $a \prec^+ b$. Moreover, since $Q\times\{1\}$ is discrete, we have that if $Q\times\{1\}\cap [z_0,y_0)_{\prec}$ or $Q\times\{1\}\cap (y_1,z_1]_{\prec}$  is infinite, then it is a sequence, of order type $\omega$, which   converges to $y_0$ or a sequence, of order type $\omega^*$, which converges to  $y_1$, respectively. Similarly, if $a,b$ is a jump in $Y$, and $Q\times\{1\}\cap(a,b)_\prec$ is infinite, then    it has order type either $\omega$ or $\omega^*$ or $\mathbb Z$, and it accumulates to $a$ or $b$ or both, respectively.  
 
\medskip 
 
For each jump $a,b$ of $Y$ such that $(a,b)_{\prec}\cap Q\times\{1\}\ne\emptyset$, let $(q_{a,b},1)=\min((a,b)_{\prec}\cap Q\times\{1\})$ if there is a minimum, and let $(q_{a,b},1)=\max((a,b)_{\prec}\cap Q\times\{1\})$ if there is a maximum and not a minimum, and choose $(q_{a,b},1)\in (a,b)_{\prec}\cap Q\times\{1\}$ arbitrarily  otherwise.
 
 \medskip
 
 Define $\tilde{f}:Q\times\{\,1\,\} \rightarrow Y$ as follows:
 $\tilde{f}(q,1)=$

 $$
 \begin{cases}
 y_0 & \textrm{if } (q,1)\prec y_0 \ \textrm{and}\ \pi_1(y_0)\ne q \\
 y_1  & \textrm{if } (q,1)\prec y_0 \ \textrm{and}\ \pi_1(y_0)= q \\
  y_1 & \textrm{if } y_1\prec (q,1) \ \textrm{and}\ \pi_1(y_1)\ne q \\
 y_0  & \textrm{if } y_1 \prec (q,1) \ \textrm{and}\ \pi_1(y_1)= q \\
 b   &\textrm{if } (q,1)\in (a,b)_{\prec^+}, (q_{a,b},1) \ \textrm{is a minimum,}\ \textrm{and} \ \pi_1(b)\ne q \\
 a   &\textrm{if } (q,1)\in (a,b)_{\prec^+}, (q_{a,b},1) \ \textrm{is a minimum,}\  \textrm{and} \ \pi_1(b)= q \\
 a   &\textrm{if } (q,1)\in (a,b)_{\prec^+}, (q_{a,b},1) \ \textrm{is a maximum,}\ \textrm{and} \ \pi_1(b)\ne q \\
 b   &\textrm{if } (q,1)\in (a,b)_{\prec^+}, (q_{a,b},1) \ \textrm{is a maximum,}\ \textrm{and} \ \pi_1(b)= q \\
 b   &\textrm{if } (q,1)\in (a,b)_{\prec^+}, (q_{a,b},1) \ \textrm{is neither a minimum or maximum,}\ \\ {}& (q_{a,b},1)\preceq (q,1) \ \textrm{and} \ \pi_1(b)\ne q \\
 a   &\textrm{if } (q,1)\in (a,b)_{\prec^+}, (q_{a,b},1) \ \textrm{is neither a minimum or maximum,}\ \\ {}& (q_{a,b},1)\preceq (q,1) \ \textrm{and} \ \pi_1(b)= q \\
 a    &\textrm{if } (q,1)\in (a,b)_{\prec^+}, (q_{a,b},1) \ \textrm{is neither a minimum or maximum,}\ \\ {}& (q,1)\prec (q_{a,b},1) \ \textrm{and} \ \pi_1(a)\ne q \\
 b   &\textrm{if } (q,1)\in (a,b)_{\prec^+}, (q_{a,b},1) \ \textrm{is neither a minimum or maximum,}\ \\ {}& (q,1)\preceq (q_{a,b},1) \ \textrm{and} \ \pi_1(a)= q \\

 \end{cases}
 $$
 and let  $f\colon Q\rightarrow [0,1]$ be given by $f(q)=\pi_1(\tilde{f}(q,1))$. It follows directly from the definition that $f$ does not have any fixed points. 
 
 \begin{claim}
 For every $x\in X\setminus( Q\cup f''(Q))$ the set $\{q\in Q: x\in [q,f(q)]\}$ is finite.
 \end{claim}
 
\begin{proofclaim}
Suppose for a contradiction that there is an $x$ such that \\ 
 $\{q\in Q\colon  x\in [q, f(q)]\}$ is infinite.  By going to an infinite subset we can find an infinite  sequence $\{q_n\colon n\in\omega\}$ 
 so that, without loss of generality,  $q_n< x< f(q_n).$  We proceed by cases.
 
 \medskip
\textbf{Case 1:} If $\{f(q_n)\colon n\in\omega\}$ is finite, then we can find an infinite subset $\Gamma\subseteq \omega$ and $a\in Y$ such that $\tilde{f}(q_n,1)=a$ for all $n\in \Gamma$. It follows, from the construction of $\tilde{f},$ that $\{(q_n,1)\colon n\in \Gamma\}$ converges to $a$. On the other hand, $\{[(q_n,0)]\colon n\in \Gamma\}\subseteq \pi'' ([(0,0),(x,0)]_{<_\ell})$  and $a\in  \pi'' ([(x,1),(0,1)]_{<_\ell}),$ which contradicts that for all infinite $A\subseteq Q,$ $A\times\{1\}$ and $\pi''(A\times\{0\})$ have the same accumulation points.
 
\medskip
\textbf{Case 2:} If $\{f(q_n)\colon n\in\omega\}$ is infinite, then we can choose an infinite convergent  $\prec$-monotone (say, increasing)   subsequence $\{q_{n_k}\colon k\in \omega\}$ so that $(q_{n_0},0) \prec^+{f}(q_{n_1},0)\prec (q_{n_2},0)\prec \dots$ Thus, both sequences converge to a point in the set
$ \pi'' ([(x,1),(0,1)]_{<_\ell}),$ which as before contradicts the fact that $\{q_{n_k}\colon k\in \omega\}$ is contained in the closed set $ \pi'' ([(0,0),(x,0)]_{<_\ell})$.
\end{proofclaim}     

\medskip

Set $X_n:=\{x\in (X\setminus (Q\cup f''(Q))\colon |\{q: x\in (q,f(q))\}|=n\}$, and let $g\colon Q\rightarrow 2$ be the function given by $g(q)=0$  iff  $f(q)<q$.

\begin{claim}
For each $n,$ there is a $G_\delta$ set $G_n$ in $(X,\leq, \tau_g)$  so that $Q\subseteq G_n$ and $G_n\cap X_n=\emptyset$.
\end{claim}

\begin{proofclaim}
We construct the desired $G_\delta$ by recursion on $n$.
If $n=0$, then set $G_0=\bigcup_{q<f(q)} [q,f(q))\cup\bigcup_{f(q)<q}(f(q),q]$, it follows from the definition of $X_0$ that $X_0\cap G_0=\emptyset$. 

Suppose  $n\geq 1$. For each $F\in [Q]^n$, let $U_F:=\bigcap_{q\in F} (q, f(q))$ and let $G=[0,1]\setminus \bigcup_{F\in [Q]^n} U_F$. For every $q\in       \bigcup_{F\in [Q]^n} U_F \cap Q$ pick an $F$ so that $q\in U_F$. Now, choose $\alpha_q <q<\beta_q$ such that  $(\alpha_q,\beta_q)\subset U_F.$  Define $G_n=G\cup \bigcup_{q<f(q)}[q,\beta_q)\cup  \bigcup_{q>f(q)}(\alpha_q,q]$. Note that $G_n$ is a $G_\delta$ in $(X,\leq, \tau_g)$, as is a finite union of $G_\delta$'s, and  it follows from the definition of $X_n$ that $G_n\cap X_n=\emptyset$.
\end{proofclaim}

\medskip

It follows that  $Q=\bigcap_{n\in\omega} G_n\cap\bigcap_{x\in (f''(Q)\cap X)\setminus Q}([0,1]\setminus\{x\}). $ Hence, $Q$  is a  $G_\delta$ in $(X,\leq, \tau_g)$ which contradicts the assumptions about $X$.
\endproof

\bigskip

 \section{Lower bound for the hereditary interval number.}
The purpose of this section is to prove the converse inequality, i.e.,  that $\hia\leq\mu$. This will be done in a series of lemmas, but before we begin let us recall a few definitions.

\medskip

Given a partially ordered set $(P,\leq)$, we denote by $\mathcal B(P)$ the Boolean subalgebra of the power-set algebra of $P$ generated by the cones of  $b_x = \{\,y\in P \colon x \leq y\,\}$ for $x\in P$.  It is easy to see that if $P$ is a linearly order set, then $\mathcal B(P)$ is an interval algebra. Moreover, one can show that every interval algebra is isomorphic to an algebra of the form $\mathcal B(P)$ for some linear order $P$.

\medskip

The class of subinterval algebras also admit a similar representation. Recall that a $pseudotree$ is a partially ordered set $T$ so that  the set $\{\,y\in T \colon y \leq x\,\}$ is linearly ordered for all $x\in T$.  A $pseudotree$ $algebra$ is a Boolean algebra of the form $\mathcal B(T)$, where $T$ is a pseudotree. The following result  relates both classes of Boolean algebras.

\begin{thm}[see \cite{Be}, \cite{He}, \cite{KoMo}]\label{pseudotreethm} 
The class of subinterval algebras coincides with the class of pseudotree algebras.\hfill $\square$
\end{thm}

The rest of the section is devoted to the proof of the following theorem.

\begin{thm}\label{mains3} 
Any $\sigma$-centered pseudotree algebra of cardinality less than $\mu$ is isomorphic to an interval algebra.
\end{thm}

It is worth mentioning that Theorem \ref{mains3} is a  generalization of Theorem 3.1 of \cite{BeTo}, we will build on their work,  and also borrow some of its notation. However, in order to carry this generalization we have to work  with the Stone space instead of working directly with the pseudotree.

\medskip

Let $(T, \leq)$ be a rooted pseudotree with root $0_T$ of cardinality $<\mu$ and such that the Boolean algebra $\mathcal B(T )$ is $\sigma$-centered.  We will show that there is a linear order $\prec_\infty$ on the Stone space $X:=Ult(\mathcal B(T))$ which induces the Stone topology. In order to do so, we define a sequence $\sim_n$ of Boolean equivalence relations on $X$, and a sequence $\prec_n$ of linear orders on the quotient spaces $X_n:= X/\sim_n$ such  that for any $x, y \in X$ and any $n\in\omega$, if  $x\sim_{n+1} y, $  then $x\sim_n y$, the linear order $\prec_n$ induces the quotient topology on $X_n$, the maps $\pi_{n+1,n}\colon X_{n+1}\to X_n$ are continuous, monotonic surjections and $X$ is the inverse limit of the inverse system $( X_n, \pi_{n+1,n}\colon n\in\omega)$, and $\prec_\infty$ is, roughly speaking, the limit of the orders $\prec_n$.
 
\medskip 
 
Recall that the Stone space of $\mathcal B(T)$ can be viewed as  the set $P(T)$ of all downward closed chains (\emph{paths}) of $(T,\leq)$ with the topology induced from the Cantor cube $2^T$, when paths are identified with the corresponding characteristic functions.  

\medskip

Notice that since  $\mathcal B(T)$ is $\sigma$-centered then  $T$ can be decomposed into countably many chains and moreover, each $C_n$ is isomorphic to a suborder of the reals. Fix a sequence $C_n (n\in\omega)$ of $\subseteq$-maximal chains of    $T$ that cover $T$ and let $P(C_n)$ be the set of all  paths $p$ contained in $C_n$ and let $P'(C_n)$ be the set of all elements of $P(C_n)$ without a supremum in $C_n$. For any subset $B\in [P'(C_n)]^{\leq\omega}$. Define a linear order $<_{n,B}$ on $C_n\cup B$ as follows: $x<_{n,B} y$ iff either
\begin{itemize}
\item  $x,y \in C_n$ and $x<y$ or
\item $x\in C_n, y\in B$ and $x\in y$ or 
\item $x\in B, y\in C_n$  and $y\notin x$ or 
\item $x,y \in B$ and $x\subsetneq y$.
 \end{itemize}

\begin{lem}\label{lemma gd}  For every $n\in\omega$ and every countable subset $B$ of the set $P'(C_n)$, there exists a  function $f\colon B\rightarrow 2$ such that    $B$ is a $G_\delta$-subset of the generalized order space $(B \cup C_n, <_{n,B},\tau_f)$.
\end{lem}
\proof Notice that since $C_n$ is isomorphic to a suborder of the reals $\mathbb R$, then $B\cup C_n$ is order-separable. The rest follows from  the definition of $\mu$.\endproof 

For $p \in P (T )$, set $T_p = \{\,y \in T \colon x < y\ \rm{for\ all}\  x \in p\,\}$. Define an equivalence relation $\equiv_p$ on $T_p$  by $x\equiv_p y$ iff there exists $z\in T_p$ such that $z\leq x$ and $z\leq y$. This lets us define the set $BP(T) = \{\,p \in P(T) \colon |T_p/\equiv_p| \geq 2\,\}$ of branching paths of $T$.

\begin{lem}[\cite{BeTo}]
$BP(T)$ is countable.\hfill $\square$
\end{lem}

Let $B_0=BP(T)\cap P(C_0)$ and let $B'''_0 =  BP(T)\cap P'(C_0)$. It will be important later on that we preserve the root  in our construction, so  if  the root is a branching point which does not have an immediate successor in $C_0$, then we  also add $0_T$ to $B_0'''$. Since $B'''_0$ is a countable subset of $P'( C_0)$ we can find,  by Lemma \ref{lemma gd},  a function $f\colon B'''_0\rightarrow 2$ (we also require that $f(0_T)=1$ in case $0_T\in B_0'''$) and  a decreasing sequence $U_n$ of open subsets of $(B'''_0 \cup C_n, <_{n,B'''_0},\tau_f)$ such that $B'''_0 =\bigcap_{n\in\omega} U_n$. 
For $p \in B_0$, let $\mathcal K_p = (T_p/\equiv_p)\setminus \{[C_0 \setminus p]\}$, where $[C_0 \setminus p]$ denotes the $\equiv_p$-equivalence class that includes $C_0 \setminus p$. Thus, $\mathcal K_p$ is the set of all classes of the quotient $T_p/\equiv p$ with the exception of the class where the points of $C_0 \setminus p$ belong. By the definition of $B_0$, we have that $\mathcal K_p\ne\emptyset$ for all $p \in B_0$.

\medskip

 Fix an enumeration $\{p_\ell\colon \ell <\ell(B_0)\}$, without repetitions, of $B_0$.  For each $\ell<\ell(B_0)$, let $\{ K_{i,\ell} \colon i<i(\ell)\}$ be an enumeration, without repetitions, of the set $\mathcal K_{p_\ell}$ and for each $K_i \in \mathcal K_{i,\ell}$,    let $x_{0,i,\ell}=\min(K_{i,\ell})$ if it exists,  and otherwise fix a strictly decreasing sequence   $ \{x_{n,i,\ell}\colon n< n(i,\ell)\}$ of elements of $K_{i,\ell} $ coinitial on $K_{i,\ell} $. We will also use the notation $x_{n,K,p}$ whenever $p\in B_0$ and $K\in \mathcal K_p$.
 
\medskip 
 
For any $p\in P(T)\setminus P(C_0)$, 
let $K^0_p=[p\setminus C_0]_{\equiv_{p\cap C_0}}$ and let $x^0_p=  x_{k, K^0_p,p\cap C_0}$, where $k=\min\{n\in\omega \colon x_{n,K^0_p,p\cap C_0} \in p\setminus C_0 \}$. Define an equivalence relation $\sim_0$ on $P(T)$ as follows:
for any $p,q \in P(T)$ we say that $p\sim_0 q$ iff either $p,q\in P(C_0)$ and $p=q$ or $p, q\notin P(C_0)$ and $x^0_p=x^0_q$.

\begin{lem}
The relation $\sim_0$ is a Boolean equivalence relation.
\end{lem}
\proof
It is straightforward  to verify that $\sim_0$ is an equivalence relation. Let us show that $\sim_0$ is Boolean.
Fix   two nonequivalent elements $p, q\in P(T)$. We need to find a clopen set $b$ which is $\sim_0$-saturated such that $p\in b$ and $q\notin b$. The proof proceed by cases:

\medskip
\textbf{Case 1:} First supose the case  $p\cap C_0\ne q\cap C_0$.  Choose $t\in q\bigtriangleup p$ arbitrary and  consider the basic clopen set $N_t:=\{r\in P(T)\colon t\in p\}$. It suffices to show that $N_t$ is $\sim_0$-saturated. Let $r\in N_t$ and $r'\in P(T)$ so that $r\sim_0 r'$. Notice that $t\in r\cap C_0=r'\cap C_0$, hence $r'\in N_t$.

\medskip
\textbf{Case 2:} Suppose  $p\cap C_0= q\cap C_0.$ Notice that  at least one of the sets $p\setminus C_0$, $q\setminus C_0$ is non-empty.   Let us assume, without loss of generality, that  that $p\setminus C_0$ is nonempty.  

\medskip
\textbf{\quad Case 2.a:} If $q\setminus C_0=\emptyset$, then consider the basic clopen set $N_p:=\{r\in P(T)\colon x^0_p\in r\}$. It is sufficient to show that $N_p$ is $\sim_0$-saturated. For any given $r\in N_p$ and $r'\in P(T)$ so that $r\sim_0 r'$. Observe that $x^0_p\leq x^0_r=x^0_{r'}$. Thus, $x^0_p\in r'$.

\medskip
\textbf{\quad Case 2.b:} If $q\setminus C_0\ne\emptyset$, then choose $x=\max\{x_p^0,x_q^0\}$ if they are comparable and arbitrarily otherwise. It follows that the basic clopen set  $N_x:=\{r\in P(T)\colon x\in r\}$ separates $p$ from $q$ and it is $\sim_0$-saturated, as shown, in the previous case.
\endproof


Let $B_0'$ be the collection of all  branching paths that have an immediate succesor in $C_0$, i.e., there exists an element $t_p\in C_0$ such that $p=\{x\in T\colon x< t_p\}$. Let $B_0''$ be the collection of all paths in $B_0\setminus (B_0'\cup\{0_T\})$ that have a maximum in $C_0$, and recall that $B_0'''$  was defined above.. 

\begin{lem}\label{lemma phi} There is a function $\varphi\colon B_0'\sqcup B_0''\sqcup B_0'''\times\omega\to C_0$ such that:
\begin{enumerate}
\item $p=[0_T,\varphi(p))$ if $p\in B_0'$;
\item $p=[0_T,\varphi(p)]$ if $p\in B_0''$;
\item $\varphi\restriction_{B_0'''\times\omega}$ is one-to-one and
\item for all $x\in C_0$, $\{\, (p,n)\in B_0'''\times\omega\colon x\in [p, \varphi (p)]\,\}$ is finite. 
\end{enumerate}
\end{lem} 

\proof First notice that clauses $(1), (2),$ implicitly define $\varphi\restriction_{B_0'\sqcup B_0''}$. Thus, it suffices to define  $\varphi\restriction_{B_0'''\times\omega}$. Fix an enumeration without repetitions $\{\,(p_k,n_k)\colon\in \omega\,\}$ of $B_0'''\times\omega$. Now for each  $k\in\omega,$ if $f(p_k)=0 $, then  choose  $\varphi(p_k,n_k)<_{0,B_0'''}p_k$ such that $[\varphi(p_k,n_k), p_k]\subseteq U_k$ and $\varphi(p_k,n_k)> \max\{\, \varphi(p_j,n_j)\colon j<k, \varphi((p_j,n_j)<_{0,B_0'''} p_k\,\}$ and in case $f(p_k)=1 $ choose  $\varphi(p_k,n_k)>_{0,B_0'''}p_k$ such that $[\varphi(p_k,n_k), p_k]\subseteq U_k$ and $\varphi(p_k,n_k)< \min\{\, \varphi(p_j,n_j)\colon j<k, \varphi((p_j,n_j)<_{0,B_0'''} p_k\,\}$. Observe that we can choose such an element as the paths in $B_0'''$ are Dedekind cuts of $C_0$. Notice that condition $(3)$ holds by construction. So we need only to check that the condition (4) is satisfied. In order to verify it, fix $x\in C_0$ and choose $N\in\omega$ so that $x\notin U_k$ for $k\geq N$.  We show that  $x\notin [p_k, \varphi (p_k,n_k)]$ for any $k\geq N$. We now proceed by cases:

\medskip
\textbf{Case a:} If $x<_{0,B_0'''} p_k$  and $f(p_k)=0$, then      the result follows as $x\notin U_k$ and $[\varphi(p_k,n_k), p_k]\subseteq U_k$.

\medskip
\textbf{Case b:} If $x<_{0,B_0'''} p_k$  and $f(p_k)=1$, then       the result follows as $p_k<_{0,B_0'''} \varphi(p_k,n_k)$.

\medskip
\textbf{Case c:} If $x>_{0,B_0'''} p_k$  and $f(p_k)=0$, then       the result follows as $p_k>_{0,B_0'''} \varphi(p_k,n_k).$

\medskip
\textbf{Case d:} If $x>_{0,B_0'''} p_k$  and $f(p_k)=1$, then       the result follows as $x\notin U_k$ and $[\varphi(p_k,n_k), p_k]\subseteq U_k$.
\endproof

\medskip

\begin{lem}
The Boolean space $X_0:=P(T)/\sim_0$ is orderable. 
\end{lem}
\proof
First note that $X$ is homeomorphic to $P(T_0)$, where $T_0:=C_0\cup \{x_{n,i,\ell}\colon n<n(j,\ell), i<i(\ell), \ell<\ell(B_0)\}$. 

\medskip

Now we define a linear order $(L, \prec_0)$ as follows: as a set $$L:= P(C_0)\cup \bigcup\{\, \zz\times\{x\}\colon x\in \ima (\phi)\,\}$$ 
we order the elements of $P(C_0)$ by inclusion and  for each $x\in\ima(\varphi)$ we insert a copy of the integers $\zz\times\{x\} $ between the paths $[0_T,x)$ and $[0_T, x]$ with its usual ordering.  The order on $P(T_0)$ will be the induced from the  embedding  $\Phi\colon P(T_0)\to L$. 
Define $\Phi$ as follows:
For each $p \in B_0'$, let $$\Phi \colon  \{\,[0_T, x_{n,K,p}]\colon n<n(K,p), K\in \mathcal L_p'\,\}\to  \zz^-\times\{\,\varphi(p)\,\}$$ be any injection onto an interval  that contains $(-1,\varphi(p))$. 
For each $p \in B_0''$, let $$\Phi \colon  \{\,[0_T, x_{n,K,p}]\colon n<n(K,p), K\in \mathcal L_p\,\}\to  \zz^+\times\{\,\varphi(p)\,\}$$ be any injection onto an interval  that contains $(1,\varphi(p))$.
Finally, for each $p\in B_0'''$, let $\Phi([0_T,x_{n,K,p}]):=(0,\varphi(p))$.
We claim that the order topology on $X_0$ coincide with the quotient topology.  First observe that, since we preserve the root, $\{\,0_T\,\}, C_0$  are the minimum, and maximum of $P(T_0)$, respectively.  Notice that $(X_0,\prec_0)$ is a complete linear order and thus, the order topology is compact Haudorff, since both topologies are compact Hausdorff, they   are $\subseteq$-minimal (among all Hausdorff topologies). Hence,  it suffices to show that these topologies are comparable. We shall show that for all $x\in T_0$ the set $N_x$ are clopen in the order topology. As we are dealing with several linear orders we will use the notation $I(p,\infty):=\{\, q\in P(T_0)\colon p\prec_0 q\,\}$ and similarly the self-explained notation $I(-\infty, p)$.   First observe, that if $x\in T_0\setminus C_0$, then $N_x$ is a finite discrete set in both topologies, so we are done in this case. We now proceed by cases:

\medskip
\textbf{Case 1a:} If $[0_T,x)\in B_0'$ and $[0_T,x]\in B_0''$, then 
$$N_x= (I([0_T,\Phi^{-1}(-1,x)],\infty)\cup F)\setminus G,$$ 
and 
$$N_x^c=(I(-\infty, [0_T,\Phi^{-1}(1,x)])\cup G)\setminus F,$$
where 
\begin{align*}
F=\{\, [0_T, t] \colon & t\in T_0\setminus C_0, [0_T,t]\cap C_0\in B_0''', \\{}& \varphi([0_T,t]\cap C_0)<_{0,B_0'''} x<_{0,B_0'''} [0_T,t]\cap C_0 \,\}
\end{align*}
and
\begin{align*}  
G=\{\, [0_T, t] \colon & t\in T_0\setminus C_0, [0_T,t]\cap C_0\in B_0''', \\{}& [0_T, t]\cap C_0<_{0,B_0'''} x\leq_{0,B_0'''}\varphi( [0_T,t]\cap C_0) \,\}.
\end{align*}
By Lemma \ref{lemma phi}, $F,G$ are finite sets.  Hence, $N_x$ is clopen in the order topology.

\medskip
\textbf{Case 1b:} If $[0_T,x)\in B_0'$ and $[0_T,x]\notin B_0''$, then 
$$N_x= (I([0_T,\Phi^{-1}(-1,x)],\infty)\cup F)\setminus G,$$ 
and 
$$N_x^c=(I(-\infty, [0_T,x])\cup G)\setminus F,$$
where 
\begin{align*}
F=\{\, [0_T, t] \colon & t\in T_0\setminus C_0, [0_T,t]\cap C_0\in B_0''', \\{}& \varphi([0_T,t]\cap C_0)<_{0,B_0'''} x<_{0,B_0'''} [0_T,t]\cap C_0 \,\}
\end{align*}
 and  
 \begin{align*}
 G=\{\, [0_T, t] \colon & t\in T_0\setminus C_0, [0_T,t]\cap C_0\in B_0''', \\{}& [0_T, t]\cap C_0<_{0,B_0'''} x\leq_{0,B_0'''}\varphi( [0_T,t]\cap C_0) \,\}.
 \end{align*}
 By Lemma \ref{lemma phi}, $F,G$ are finite sets.  Hence, $N_x$ is clopen in the order topology.

\medskip
\textbf{Case 2a:} If $[0_T,x)\notin B_0'$ and $[0_T,x]\in B_0''$, then 
$$N_x= (I([0_T,x),\infty)\cup F)\setminus G,$$ 
and 
$$N_x^c=(I(-\infty, [0_T,\Phi^{-1}(1,x)])\cup G)\setminus F,$$
where 
\begin{align*}
F=\{\, [0_T, t] \colon & t\in T_0\setminus C_0, [0_T,t]\cap C_0\in B_0''', \\{}& \varphi([0_T,t]\cap C_0)<_{0,B_0'''} x<_{0,B_0'''} [0_T,t]\cap C_0 \,\}
\end{align*}
 and  
 \begin{align*}
 G=\{\, [0_T, t] \colon & t\in T_0\setminus C_0, [0_T,t]\cap C_0\in B_0''', \\{}& [0_T, t]\cap C_0<_{0,B_0'''} x\leq_{0,B_0'''}\varphi( [0_T,t]\cap C_0) \,\}.
 \end{align*}
  By Lemma \ref{lemma phi}, $F,G$ are finite sets.  Hence, $N_x$ is clopen in the order topology.

\medskip
\textbf{Case 2b:} If $[0_T,x)\notin B_0'$ and $[0_T,x]\notin B_0''$, then 
$$N_x= (I([0_T,x),\infty)\cup F)\setminus G,$$ 
and 
$$N_x^c=(I(-\infty, [0_T,x])\cup G)\setminus F,$$
where 
\begin{align*}
F=\{\, [0_T, t] \colon & t\in T_0\setminus C_0, [0_T,t]\cap C_0\in B_0''', \\{}& \varphi([0_T,t]\cap C_0)<_{0,B_0'''} x<_{0,B_0'''} [0_T,t]\cap C_0 \,\}
\end{align*}
 and  
\begin{align*} 
G=\{\, [0_T, t] \colon & t\in T_0\setminus C_0, [0_T,t]\cap C_0\in B_0''', \\{}& [0_T, t]\cap C_0<_{0,B_0'''} x\leq_{0,B_0'''}\varphi( [0_T,t]\cap C_0) \,\}.
\end{align*}
By Lemma \ref{lemma phi}, $F,G$ are finite sets.  Hence, $N_x$ is clopen in the order topology.
It follows that both topologies coincide. \endproof

   Now we define a finer equivalence relation $\sim_1$ using the chain $C_1$.  Let $p=C_0\cap C_1$ and let $K$ be the class of $T_p/\equiv_p$ which contains $C_1$.
   If $K$ has a minimum, then $K=[x_{0,K,p}]_0$ is a clopen set in $X$ and consider the rooted pseudotree $T_1=\{\, t\in T\colon t\geq x_{0,K,p}\,\}$. We identify the paths of $P(T_1)$ with the elements of $N_{x_{0,K,p}}$ via the map $p\mapsto [0_t,x_{0,K,p})\cup p$. Define an equivalence relation $\sim_1$ on the paths  $P(T_1)$, as before,  using  the path $C_1\setminus [0_T,x_{0,K,p})$  and we also obtain a linear order $\prec_1$.
 When $K$ has no minimal element, we apply the previous process, separately, to the clopen sets $[x_{n,K,p}]_0$ using the rooted pseudotrees $T_1^n:=\{\, t\in T\colon t\geq x_{n,K,p}, \ t\perp   x_{n-1,K,p}\,\}$.
 
\medskip 
 
 Proceeding in this way, we construct a sequence of finer equivalence relations   $\sim_n$ and linear orders $\prec_n$ on $X_n:=X/\sim_n$ such that for all $n\in\omega$,
 \begin{itemize}
 \item for all $p, q\in X$, $p\sim_n q$ implies $p\sim_{n+1} q$,
 \item the order topology induced from $\prec_n$ coincide with the quotient topology on $X_n$,
 \item the natural map $\pi_{n+1,n}\colon X_{n+1}\to X_n$, $[x]_{n+1}\mapsto [x]_n$  is  a continuous non-decreasing surjection.
  \end{itemize}
 Let $X_\infty$ be the inverse limit of the inverse system $(   \pi_{n+1,n}\colon X_{n+1}\to X_n; n\in\omega)$ and define a linear order $\preceq_\infty$ on $X_\infty$ as follows
 $x\preceq_\infty y$ iff $\pi_n(x)\preceq_n \pi_n(y)$ for all $n\in \omega$.

\medskip

 We are now ready to finish the  proof of Theorem \ref{mains3}.
   
   \begin{lem}
   $X$ is homeomorphic to $X_\infty$  and the order $\prec_\infty$ induces the Cantor topology on $X$.
    \end{lem}
\proof First we verify that $X$ is homeomorphic to $X_\infty$. Define $h\colon X\to X_\infty$ by sending $x\to h(p)(n)=[p]_n$. It is clear that $h$ is continuous and since $X$ is compact, it suffices to show that $h$ is a bijection. Let us first verify that $h$ is one-to-one. Fix $p\ne q\in X$ and pick an element $t\in p\bigtriangleup q$. Find $n\in\omega$ such that $t\in C_n$. It follows from the definition of $\sim_n$ that $p \not\sim_n q$. Now let us  check that $h$ is onto let $x\in X_\infty$. Pick $p_n\in X$ so that $x(n)=[p_n]_n$. Notice that $[p_{n+1}]_{n+1}\subseteq [p_n]_n$ since the equivalence classes are closed and $X$ is compact it follows that there is a $p\in \bigcap_n [p_n]_n$. Thus, $h(p)=x$ as required.

\medskip

We are left to show that the order topology coincide with the Cantor topology. In order to do so, we  first prove that  $(X,\prec_\infty)$ is a complete linear order.  Let $A\subseteq X$ be given. Consider the set  $A_n;=\{\, [x]_n\colon x\in A\,\}$, and  choose  $a_n\in X$ so that $[a_n]_n=\sup A_n$. Observe that for any $n\in\omega, [a_{n+1}]_{n=1}\subseteq [a_n]_n$. It follows that  $\{\,a\,\}:=\bigcap_{n\in\omega} [a_n]_n$ is equal to the supremum of $A$.

\medskip

Let us now show that both topologies coincide,  by using the same argument  as above,  it is sufficient to show that the initial (final) segments $I(-\infty, p)  (I(p,\infty)) $ are open in the Cantor  topology.
Notice that $  I(-\infty, p)=\bigcup_n \pi_n^{-1}(I(-\infty, [p]_n)$. Thus,  $I(-\infty, p) $ is an open set. Analogously, the final segments are also open. This concludes the proof of the Lemma.   \endproof 

%
\section{Comparison with classical cardinal invariants.}
In this section we shall compare the cardinal $\hia$ with the cardinals in the Ciho\'n's diagram and we  also show that it is consistent that $\b<\hia$.

\begin{lem}
$\b \leq \hia\leq \min\{\, \mathfrak d, non(\mathcal M)\,\}$. 
\end{lem}
\proof
As mentioned before the  first inequality follows from \cite{BeTo}. The inequality $\hia\leq non(\mathcal M)$ follows from the trivial observation that for any $f\in 2^{\qq}$ any $G_\delta$-set in the $(X\cup \qq, \tau_f,\leq)$ which contains $\qq$ is comeager.
Let  us now verify that $\hia\leq \d$. Fix a $\leq$-dominating family $\mathcal F$ of functions $g\colon \qq \to \omega\setminus\{\, 0\,\}$ of cardinality $\d$, we may also assume that for every $g\in\mathcal F$ the set  $|g^{-1}([0,n])|$ is finite  for all $n\in\omega$.   For each $g\in\mathcal F$,  recursively construct   a sequence of rational numbers $( q_s\in\qq\colon s\in 2^{<\omega})$ and a Cantor's scheme  of closed intervals $(I_s\colon s\in 2^{<\omega})$ such that
$I_s=[-\epsilon_s+q_s,q_s+\epsilon_s]\subseteq (-\frac{1}{g(q_s)}+q_s, q_s+\frac{1}{g(q_s)})$,  $I_{s^\smallfrown0}\subseteq (-\epsilon_s+q_s,q_s)$, and $I_{s\smallfrown 1}\subseteq ( q_s,q_s+\epsilon_s)$.  Choose $X\subseteq \rr\setminus \qq$ of cardinality $\aleph_1$ so that for each $x\in X$ there is a $s_x\in 2^\omega$ so that $\{\,x\,\}=\bigcap_{n\in\omega} I_{s_x\restriction n}$. Let us write $X_g$ as  $X_g:=\{\,x_\alpha\colon \alpha\in\omega_1\,\}$. For each $\alpha\in\omega_1$, set $a_\alpha:=\{q\in \qq\colon 0<x_\alpha-q<\frac{1}{g(q)}\}$, and $b_\alpha:=\{q\in\qq\colon 0<q-x_\alpha<\frac{1}{g(q)}\}$.  It follows from the properties of the  construction  that $a_\alpha, b_\alpha$ are both infinite and disjoint. Moreover, the family $(a_\alpha, b_\alpha\colon \alpha\in\omega_1)$ forms a Luzin gap (see \cite{Stevo1}). In order to verify this,   let $\alpha\ne\beta\in\omega_1$ be given, and let $n=\Delta(s_x,s_y)$. We have that    either $q_{s_x\restriction n}\in a_\alpha\cap b_\beta$ or $q_{s_x\restriction n}\in a_\beta\cap b_\alpha.$
Now set $X:=\bigcup_{g\in\mathcal F} X_g$. Notice that  clearly $|X|\leq \d$.

\begin{claim}
For any $f\in 2^{\qq}$, the rational numbers $\qq$ are not $G_\delta$ in the generalized order space  $(X\cup \qq, \tau_f,\leq).$
\end{claim}

\begin{proofclaim}
Suppose for a contradiction that there is a function $f\in 2^{\qq}$, and a $G_\delta$-set $G$ which contains $\qq$ so that $G\cap(X\cup\qq)=\qq$. Let $c=f^{-1}(0)$. Using the fact that $\mathcal F$ is dominating we can find $g\in \mathcal F$ such that    $$\tilde G:=\bigcap_{n\in\omega}(\bigcup_{q\in c} (-\frac{1}{\max\{\,n,g(q)\,\}}+q,q]\cup    \bigcup_{q\notin c} [q,q+\frac{1}{\max\{\,n,g(q)\,\}})\subseteq G.$$
Notice that if $|c\cap a_\alpha|=\omega,$ then $x_\alpha\in  \bigcap_{n\in\omega}(\bigcup_{q\in c} (-\frac{1}{\max\{\,n,g(q)\,\}}+q,q])$, and analogously if $|(\qq\setminus c)\cap b_\alpha|=\omega,$  then $x_\alpha\in \bigcap_{n\in\omega}(  \bigcup_{q\notin c} [q,q+\frac{1}{\max\{\,n,g(q)\,\}})$. Thus, for any $\alpha\in\omega_1$ we have that $a_\alpha\cap c=^*\emptyset$ and $b_\alpha\subseteq_* c$, i.e., $c$ separates the Luzin gap $(a_\alpha, b_\alpha\colon \alpha\in\omega_1)$ which is impossible.
\end{proofclaim}
\endproof

\medskip

Before proceeding any further we will make a straightforward reformulation of the cardinal invariant $\hia$ in terms of the Cantor set.  Recall that given any $A\subseteq 2^{<\omega}$
the $G_\delta$-closure $\pi''(A):=\{\, x\in 2^\omega\colon \exists^\infty n \ x\restriction n\in A\,\}$.   Let $Q=(0,1)\cap \qq$, and $Y:=[0,1]\times\{\,0\,\}\cup Q\times\{\,1\,\}$. Observe that 
for any $f\in 2^Q$, the generalized order space $(I,\leq, \tau_f)$ can be identified with the set $$Z:=Y\setminus (\{\, (q,1)\in Q\times\{\,1\,\}\colon f(q)=0\,\}\cup   \{\, (q,0)\in Q\times\{\,0\,\}\colon f(q)=1\,\}),$$
with the subspace topology in $Y$. Also notice that $(Y,\leq)$ is isomorphic to $(2^\omega, <_{lex})$, as both as topological spaces, and as linear orders. Finally, recall that a set $A\subseteq 2^{<\omega}$ is adequate if there is a $\sigma:2^{<\omega}\rightarrow 2$ such that for every $s\in 2^{<\omega}$ there are infinitely many $n$ such that $s^\frown \sigma(s)^n\in A$.

\medskip 

Therefore, in light of the above discussion we obtain that:

\begin{thm}
$\hia=\min\{ |X| : X\subseteq 2^\omega\setminus Q, \ X\cap \pi''(A)\ne\emptyset \ {\rm for\  all  \ adequate}\   A\}.$

\end{thm} 

Armed with  the previous theorem we shall obtain a more precise upper bound of $\hia$ using a variation of $\b$. In order to do so, we recall a few definitions.

\begin{defn}
An \emph{interval partition} is a partition of $\omega$ into (infintely many) finite intervals $I_n (n\in\omega)$. We always assume that the intervals are numbered in the natural order, so that, if $i_n$ denotes the left endpoint of $I_n$ then $i_0=0$ and $I_n=[i_n,i_{n=1}).$ We say that the interval partition $\{\, I_n\colon n\in\omega\,\}$ \emph{dominates} another interval partition $\{\, J_n\colon n\in\omega\,\}$ if $\forall^\infty n\, \exists k\, (J_k\subseteq I_n).$
\end{defn}

The first part of the following theorem is a well-known  reformulation of $\b$ (see \cite{Blass}). The second part was pointed out to us by Osvaldo G\'uzman-Gonz\'alez, and it is probably folklore. However, since we were unable to find a reference we give a proof for convenience  of the reader. 
\noindent\begin{thm}

 $\b$ is the smallest cardinality of any family of interval partitions not all dominated by a single partition. 
Equivalently, $\b$ is the smallest cardinality of a family $\mathcal F$ of interval partitions such that for any interval partition 
 $\{\, J_n\colon n\in\omega\,\}$ there is an interval partition 
 $\{\, I_n\colon n\in\omega\,\}\in \mathcal F$ so that 
$\exists^\infty n\, \exists k\, (J_k\subseteq I_n)$.  
\end{thm}
\proof
Fix a family $\mathcal F$ of interval partitions not all dominated by a single partition of cardinality $\b$.  For each   $\mathcal  I=\{\,I_n\colon n\in\omega\,\}\in \mathcal F$, let  
$\mathcal I^e=\{\, I_{2n}\cup I_{2n+1}\colon n\in\omega\,\}$ and $\mathcal I^o=\{\, I_0\,\}\cup\{I_{2n+1}\cup I_{2n+2}\colon n\in\omega\,\}$. Set $\mathcal G=\{\, \mathcal I^e\colon \mathcal I\in \mathcal F\,\}\cup  \{\, \mathcal I^o\colon \mathcal I\in \mathcal F\,\}$. We claim that $\mathcal G$ is as required. Let $\{\, J_n\colon n\in\omega\,\}$ be any given interval partition. By the first part, we can find $\mathcal I\in \mathcal F$ so that $\mathcal I$ is not dominated  by $\{\, J_n\colon n\in\omega\,\}$. For any given $J_n$ that does not contain any interval of $\mathcal I$, let $k_n$ be the index such that $j_n\in I_k$ if $j_{k+1}\geq i_{n+1}$, then $J_n\subseteq I_{k_n}$. On the other hand, if $j_{k_n+1}< i_{n+1}$, then $j_{k_n+2}>i_{n+1}$, as $I_{k_n+1}\not\subseteq J_n$. Thus, in either case $J_n\subseteq I_{k_n}\cup I_{k_n+1}$. It follows that $\mathcal I^e$ or $\mathcal I^o$ work, depending on weather there are infinitely many indexes $k_n$ which are even or odd. \endproof
We consider the following variant of $\b$.
\begin{defn}\label{b2}
Let $\b_2$ be the smallest cardinality of a family $\mathcal F$ of interval partitions such that for any interval partition  $\{\, J_n\colon n\in\omega\,\}$ there is an interval partition  $\{\, I_n\colon n\in\omega\,\}\in \mathcal F$ so that 
$\exists^\infty n\, \exists k_1, k_2\, (J_{k_1}\subseteq I_n\,\wedge \,J_{k_2}\subseteq I_{n+1}\, )$.  

\end{defn}

\begin{prop}
$\hia\leq \b_2$.
\end{prop}
\proof
Fix a family $\mathcal F$ of interval partitions of cardinality $\b_2$ satisfying Definition \ref{b2}. For each   $\mathcal I:=\{\, I_n\colon n\in\omega\,\}\in \mathcal F$, let $x_{\mathcal I}\in 2^\omega$ be defined as follows $x_{\mathcal I}(k)=0$ iff $k\in I_{2n}$ for some $n\in\omega$.   Let $X=\{\, x_\mathcal I\colon \mathcal I\in\mathcal F\,\}$.  It suffices to show that  $X\cap \pi''(A)\ne\emptyset $ for  all  adequate  sets  $A$.  Fix an adequate set $A$ and $\sigma\colon 2^{<\omega}\rightarrow 2$ such that for every $s\in 2^{<\omega}$ there are infinitely many $n$ such that $s^\frown \sigma(s)^n\in A$. Recursively construct an interval partition $\{\, J_n\colon n\in\omega\,\}$ as follows. Let $j_0=0$, and suppose that $j_n$ has been constructed. For any $s\in 2^{ j_n}$, let $\ell_s$ be so that $s\frown 0^{\ell_s}\in A$ if there is such an $\ell$, otherwise set $\ell_s=0$ and let $k_s$ be so that $s\frown 1^{k_s}\in A$ if there is such a $k$, otherwise set $k_s=0$.  Define $m_s=\max\{\, \ell_s, m_s\,\}$. Finally let $j_{n+1}=1+j_n=\max\{\, m_s\colon s\in 2^{ j_n}\,\}$. We claim that for any $n$ so that $\exists k_1, k_2\, (J_{k_1}\subseteq I_n\,\wedge \,J_{k_2}\subseteq I_{n+1}\, )$, we can find   a $k\geq i_n$ so that $x_{\mathcal I}\restriction k\in A$.  Fix such an $n$, and let $k_1, k_2$ minimal such that $J_{k_1}\subseteq I_n$ and $J_{k_2}\subseteq I_{n+1}$. Let  $s=x_{\mathcal I}\restriction i_{k_1}$. It follows from the construction of $x_{\mathcal I}$ that if $s\frown \epsilon^{m_s} $ extends an element of $A$ and $x_{\mathcal I}(j_{k_1})=\epsilon $, then $x_{\mathcal I}$ extends an element of $A$ of length at least $j_{k_1}+1$. On the other hand, if this is not the case then all the extensions of $s$ prefer the other color and hence,   $x_{\mathcal I}$ extends an element of $A$ of length at least $j_{k_2}+1$. \endproof
\begin{defn}
Given $s\in 2^{<\omega},$ and $(m,\epsilon)\in\omega\times 2,$ let $s\ast (m,\epsilon)$ denote the sequence $s^\frown \epsilon^m$.
\end{defn}

The rest of the section is devoted to the proof of the following theorem.
\begin{thm} Given any uncountable regular cardinal $\kappa$, it is consistent that $\mu=\kappa$ and $\mathfrak b=\aleph_1$. 

\end{thm}
\proof
Let $Q$ denote the collection of   all elements of the Cantor's set which are eventually constant. Let $\mathbb P$ denote the forcing notion consistent of elements of the form $p=\langle \sigma_p, F_p\rangle$ where  $\sigma_p: 2^{\leq n_p}\rightarrow \omega\times 2$ for some $n_p\in\omega$  and $F_p\in [2^{\leq \omega}\setminus Q]^{<\omega}$ and for any $f,g\in F_p$ there is an $i\in dom(f)\cap dom(g)$ such that $f(i)\ne g(i)$. It is worth pointing out that, we allow finite functions in $F_p$ to make the forthcoming arguments work through. 

\medskip

We  order $\mathbb P$ as follows: $p\leq q$ if $\sigma_p \sqsupseteq \sigma_q$, for all $f\in F_q$ there is a $g\in F_p$ so that $f\sqsubseteq g$, and for every $f\in F_q$ and $k\in\omega$, if $f\restriction k\in dom(\sigma_p)\setminus dom(\sigma_q)$ then $f\not \sqsupseteq  f\restriction k\ast \sigma_p(f\restriction k)$.

\medskip

Clearly $\mathbb P$ is a $\sigma$-centered poset  since  $\langle \sigma, F\rangle$ and $\langle \sigma, G\rangle$ are compatible with common extension $\langle \sigma, H\rangle$, where $H$ consists of the end nodes of the tree $(F\cup G, \sqsubseteq)$.

\medskip

If $\mathcal G$ is a $\mathbb P$-generic filter over $V$, and $\sigma_{\mathcal G} :=\bigcup \{\,\sigma \colon \exists F (\langle \sigma, F\rangle\in \mathbb P\,\}$ and $A=\{\,s^\frown \epsilon^n\colon s\in 2^{<\omega}, n\geq m,\ \rm{for}\ (m,\epsilon)=\sigma_{\mathcal G}(s)\,\}$ then $V[\mathcal G]\models (2^\omega\setminus Q)\cap V\cap \pi''(A)=\emptyset.$
 thus we can make $\mu=\kappa$ by iterating, over a model of CH, the poset $\mathbb P$, $\kappa$ many times with finite support.    So it suffices to show that in the final model, the ground model reals are still unbounded. 
 
\medskip
 
In order to do this, we use some notions and techniques from \cite{Br95}. Given a poset $\mathbb P$, a function $h:\mathbb P\rightarrow \omega$ is called a   $height\ function$ iff $p\leq q$ implies $h(p)\geq h(q)$ for any $p,q\in \mathbb P$. A pair $(\mathbb P,h)$ is $soft$ iff $\mathbb P$ is a poset, $h$ is a height function and the following three conditions hold:
\begin{enumerate}
\item[(I)] If $p_n (n\in\omega)$ is a decreasing sequence and $\exists m\in\omega\ \forall n\in\omega (h(p_n)\leq m),$ then there is a $p\in \mathbb P$ so that $p\leq p_n$ for all $n\in\omega$.
\item[(II)] Given $p,q\in\mathbb P$ and $m\in\omega$ there is $\{q_i:i\in\ell\}\subseteq \mathbb P$ so that 
 \begin{enumerate}
  \item [(i)] for all $i\in\ell$ ($q_i\leq q$ and $q_i\perp p$),
  \item [(ii)] whenever $q'\leq q$ is incompatible with $p$ and $h(q')\leq m$ then there is $i\in\ell$ such that $q'\leq q_i$.
 \end{enumerate}
\item[(III)] if $p,q$ are compatible, there is an $r\leq p,q$   with $h(r)\leq h(p)+h(q)$.
\end{enumerate}

  \begin{lem}[\cite{Br95}] Suppose $\mathbb P$ is a ccc poset, $h$ is a height function on $\mathbb P$, and $(\mathbb P,h)$ is soft then any unbounded family of functions in $\omega^\omega\cap V$ is still unbounded in $V[G]$, where $G$ is $\mathbb P$-generic over $V$. 
\end{lem}

Thus, we are left to show. 

 \begin{lem}
 $\mathbb P$ is a soft poset.
 \end{lem}
 \proof
 Let $h:\mathbb P\rightarrow \omega$ be defined by 
 $$h(p):=\max\{|\sigma_p|, \max(ran(\pi_0\circ \sigma_p), |F_p|\}.$$ 
 
 Obviously $h$ is a height function.  Let us verify that $\mathbb P$ satisfies condition (I) of the definition of soft. Given a decreasing sequence $\langle \sigma_n, F_n\rangle (n\in\omega)$ of conditions   of bounded height, it becomes eventually constant in the first coordinate, let say $\sigma$,  and the cardinality of the set of functions is also eventually constant. For instance, assume that $F_n=\{f^n_i\i\in \ell\}$ for $n\geq k$, where the enumeration is given by the  lexicographical order. It follows that for $n\geq m\geq k,$ $f^m_i\sqsubseteq f^n_i$ for all $i\in\ell$. Set $f_i:=\bigcup_{n\geq k} f_i^n$, then $\langle \sigma, \{f_i\colon i\in \ell\}\rangle$ is a lower bound.   Condition (III) is trivial, so we are left with (II).
 
 \medskip
 
 Before going into the proof, we introduce some notation. Given $\sigma\sqsubset \tau$ and $F\in[2^{\leq \omega}\setminus Q]^{<\omega}$, we say that $\tau$ \emph{respects} $F$ if   
 for every $f\in F$ and $k\in\omega$, if $f\restriction k\in dom(\tau)\setminus dom(\sigma)$ then $f\not \sqsupseteq  f\restriction k\ast \sigma_p(f\restriction k)$.
 
 \medskip
 
 Let $p=\langle \sigma, F\rangle, q=\langle \tau, G\rangle \in\mathbb P$ be given. Note that $p$ and $q$ are compatible iff $\sigma\sqsubseteq \tau $  ( resp. $ \tau\sqsubseteq \sigma$)  and $\tau $ respects $F$ (resp.  $\sigma$ respects $G$).  Hence, $p$ and $q$ are incompatible iff either $\sigma$ is incompatible with $\tau$ or  $\sigma\sqsubseteq \tau $ (resp.  $\tau\sqsubseteq \sigma$) and $\tau$ does not respect $F$ (resp.  $\sigma$ does not respect $G$). 
 
 \medskip
 
 If $q\perp p$ or $q\leq p$, then condition (II) is trivial.  So assume that $p$ and $q$ are compatible and $q\not\leq p$. We now describe how to construct the desire finite set.
 
 \begin{enumerate}
\item[(i)]  Suppose $\sigma\sqsubseteq \tau$. Then we take all the conditions of the form $\langle \tau',G\rangle$ extending $q$ such that  $dom(\tau')\subseteq m$ and $ran(\tau')\subseteq m\times2$ such that $\tau'$ does not respect $F$.
\item[(ii)]  Assume $\tau \sqsubset \sigma$.  We take all the conditions of the form $\langle \tau',G'\rangle$ extending $q$ such that $dom(\tau')\subseteq m$ and $ran(\tau')$, $|G'|\leq m$ and either 
\begin{enumerate}
 \item[(a)] $G=G'$ and $\tau'$ is incompatible with $\sigma$ or
 \item[(b)] $\tau'\sqsubset \sigma$ and $G'=G\cup\{s\ast \tau(s)\}$  for some $s\in dom(\sigma)\setminus dom(\tau')$ or
 \item[(c)] $G'=G$ and $\sigma \sqsubset \tau'$ and $\tau'$ does not respect $F$.
\end{enumerate}
\end{enumerate}
 
We have, for each case, a finite set of conditions and it is straightforward to verify that   each condition of height $\leq m$ below $q$ and incompatible with $p$ is  below one of the described conditions.   This finishes the proof of the theorem.\endproof

\section{ ZFC results }

In this section we show that the interval algebra over a Hausdorff gap is hereditary. 

\medskip

Recall that a subspace $X$ of the reals is a $\lambda'$-\emph{set} if every countable set $D\subseteq \rr$ is a relative $G_\delta$ in $X\cup D$. Let $(a_\alpha,b_\alpha\colon \alpha\in\omega_1)$ be a Hausdorff gap. By making finite  modifications to the elements of the gap we may assume that $L:=\{\, a_\alpha\colon \alpha\in\omega_1\,\}\cup\{\,b_\alpha\colon \alpha\in\omega_1\,\}$ is dense in $2^\omega$. Thus, the order and subspace topologies on $L$ coincide. 

\medskip

It is well known (see \cite{Mi})  that $L$ is a $\lambda'$-set.  The rest of the section is devoted to the proof of the following theorem.  
\begin{thm}
The interval algebra $\mathcal B(L)$ is hereditary.
\end{thm}

The proof will be done in a series of Lemmas. Let $\mathcal A$ be a subalgebra of $\mathcal B(L)$, by Theorem \ref{pseudotreethm}, there is a rooted  pseudotree $T\subseteq \mathcal B(L)$ (ordered by inclusion) such that $\mathcal A\cong \mathcal B(T)$. Since $\mathcal A$ is $\sigma$-centered we can decompose $T$ into a countable union of chains. It follows from  the results in Section 2, specifically Lemma \ref{lemma phi}, that  it suffices to show that given any chain $C$ in $\mathcal B(L)$ and any countable set $D$ of Dedekind cuts of $C$ the set $D$ is a $G_\delta$ in $B\cup C$ with the order topology. Let us first show a weaker version of this result.
\begin{lem}\label{lemalset}
Let $X\subseteq L,$ and $D\subseteq 2^\omega$ a countable set. Then $D$ is a $G_\delta$ in $X\cup D$ with the order topology.  
\end{lem}
\proof
For each $x\in X$  define $I_X(-\infty,y):=\{\, y\in X\colon y<x\,\}$, and similarly define $I_X(y,\infty)$. Without loss  of generality, we may assume that the elements of $D$ in $X\cup D$ are not isolated, as the set of isolated points is open and can be easily handled.   Now for any $d\in D$ let $d^-:=\sup   I_X(-\infty,d)$  and $d^+:= \inf I_X(d,\infty)$ (we take both the supremum  and the infimum  inside $2^\omega$). Set $\tilde D=\{\, d^-\colon d\in D\,\}\cup \{\, d^+\colon d\in D\,\}$. Since $L$ is a $\lambda'$-set, we can find a sequence of open sets $U_n (n\in\omega)$ so that $\bigcap_{n\in\omega} U_n\cap (L\cup \tilde D)=\tilde D$. For each $n\in\omega, d^-, d^+\in \tilde D$, choose $x,y,z,w \in L$ so that  $x<d^-<y$, $y<d^+<w$ and $(x,y)\cup (z,w)\subseteq U_n$. Now choose $a^n_d, b^n_d\in X$ such that $x<a^n_d<d^-$ and $d^+<b^n_d<w$.  If either $d^+$ or $d^-$ are maximum or a minimum just do the construction from one side. Define $V_n=\bigcup_{d\in D}(a^n_d,b^n_d)$. It is easy to verify that $\bigcap_{n\in\omega} V_n\cap (X\cup D)=D$.    \endproof
 
 \begin{lem}
 Let $C$ be a linear order set,  $D$ a countable set of Dedekind cuts of $C$, and let $C=\bigcup_{n\in\omega} C_n$. If $D$ is  $G_\delta$ in $C_n\cup D$ with the order topology, then $D$ is    a  $G_\delta$ in $C\cup D$ with the order topology. 
 \end{lem}
 \proof Choose $G_n$ to be a $G_\delta$ in $C_n\cup D$ with the order topology which satisfies that $G_n\cap (C_n\cup D)=D$. It follows that $G=\bigcap_{n\in \omega} G_n$  is a $G_\delta$ in $C\cup D$ and it is as required.\endproof
 
 We are now ready to prove.
 
 \begin{lem}
 Let $C\subseteq \mathcal B(L)$ be a chain in the Boolean order, and let $D$ be a countable set of Dedekind cuts of $C$.  Then $D$ is a $G_\delta$ in $C\cup D$ with the order topology.
 \end{lem}
\proof
We may assume that $C$ is uncountable as otherwise the result is trivial. Since  every element  $t\in C$ can be written in the form $[a^t_0, a^t_1)\cup\dots\cup [a_{2m(t)-2}^t, a_{2m(t)-1}^t)$, where $a^t_0<a^t_1<\dots<a_{2m(t)-2}^t< a_{2m(t)-1}$. Observe that, by our choice of the gap,  the intervals $[a_2k,a_{2k+1})$ are infinite for $k<m$.  We can decompose $C$ as the countable union of the sets $C_n:=\{\, t\in C\colon n(t)=t\,\}$.  Now for each $t\in C_n$, we can find rational numbers $q^t_k, r^t_k \ (k<2n-1)$ such that 
$a^t_0<q_0^t<r^t_0<a^t_1<q^t_1<r^t_1<a^t_2<q^t_1<r^t_1<\dots<a^t_{2m-2}<q^t_{2m-2}<r^t_{2m-2}<a^t_{2m-1}$.  We call such a set a frame of $t$, (see \cite{BeBo}). For each $F\in\qq^{4m-2}$, we can further decompose each $C_n$ into  countably many pieces $C_{n,F},$ where $C_{n,F}$ is the set of all elements of $C_n$ which have $F$ as a frame.   By the previous Lemma, we may reduce it to the case $C=C_{n,F}$.  Observe that, $D$ might no longer be a Dedekind cut in $C_{n,F}$.  However, without loss of generality, we may assume that $D$ does not have isolated points.  Notice that, by construction, given any $t,s\in C, s\subseteq t$ iff $a_2k^t<a_2k^s$, and $a_{2k+1}^s<a_{2k+1}^t$ for $k< m$. Set $A_0=\{\, a\in L\colon \exists t\in C, \ a=a^t_0\,\}$. For each $d\in D$, let $A_0^-(d):=\{\, a\in L\colon \exists t\in C\ t< d, \ a=a^t_0\,\}$ and let $A_0^+(d):=\{\, a\in L\colon \exists t\in C\ t> d, \ a=a^t_0\,\}.$ Define $d_0^-:=\sup A_0^-(d),$ and  $d_0^+:=\inf A_0^+(d),$ where the supremum and infimum are taken inside $2^\omega$.  Define $\tilde D^0=\{\, d_0^-\colon d\in D,\ d_0^-\emph{is not a maximum},\}\cup \{\, d_0^+\colon d\in D, d_0^+\emph{is not a minimum}\,\}$.  It follows from Lemma \ref{lemalset} that there is a $G_\delta$  set $G^0$, say $G^0=\bigcap_{n\in\omega} U^0_n $, in $A_0\cup \tilde D$ so that $G\cap  (A_0\cup \tilde D)=\tilde D$. 

\medskip

Now for each $d\in D,$ and $n\in \omega$, define $t^-_{0,d,n}, t^+_{0,d,n}\in C$ as follows:

\medskip
\textbf{Case 1a:} If $d_0^-$ is a maximum, let  $t^-_{0,d,n}$ to be any  element $t\in A^-(d)$ such that $a_0^t=d_0^-$ (in this case the sequence is constant);

\medskip
\textbf{Case 1b:} If $d_0^-$ is  not a maximum, let  $t^-_{0,d,n}$ to be any  element $t\in A^-(d)$ such that $(a_0^t,d_0^-]\subseteq U^0_n$;

\medskip
\textbf{Case 2a:} If $d_0^+$ is a minimum, let  $t^+_{0,d,n}$ to be any  element $t\in A^+(d)$ such that $a_0^t=d_0^+$ (in this case the sequence is constant);

\medskip
\textbf{Case 2b:} If $d_0^+$ is  not a  minimum, let  $t^+_{0,d,n}$ to be any  element $t\in A^+(d)$ such that $[d_0^+,a^t_0)\subseteq U^0_n$.

\medskip

We now carry out an analogous  argument by considering  the sets $A_1=\{\, a\in L\colon \exists t\in C,\ a^t_1=a\,\}$, and for each $d\in D$, let $A^-_1(d)=\{\, a\in L\colon \exists t\in C,\ t<d,\ a=a^t_1\,\}$  in case $d_0^-$ is not a maximum, otherwise let   $A_1^-(d)=\{\,  a\in L\colon \exists t\in C,\ t<d,\ a=a^t_1,\ a^t_0=d_0^-\,\}$, and  let $A^+_1(d)=\{\, a\in L\colon \exists t\in C,\ t>d,\ a=a^t_1 \,\}$  in case $d_0^+$ is not a minimum, otherwise let   $A_1^+(d)=\{\,  a\in L\colon \exists t\in C,\ t>d,\ a=a^t_1,\ a^t_0=d_0^+\,\}$.
Proceeding in this way, we construct a sequence $t^\pm_{\ell,d,n}\  ( \ell<2m, d\in D, n\in\omega)$ of elements of $C$, suborders $A_\ell\ (\ell<2m)$ of $C $,  and   $G_\delta$'s $G^\ell (=\bigcap_{n\in\omega}U^\ell_n) \ (\ell<2m)$ such that:

\begin{itemize}
\item $(\tilde D^\ell\cup A_\ell)\cap G^\ell=\tilde D^\ell$, for $\ell<2m$;
\item $a^{t^-_{\ell,d,n}}_\ell=d^-_\ell$ if $d^-_\ell$ is maximum;
\item $(a^{t^-_{\ell,d,n}}_\ell, d^-_\ell]\subseteq U^\ell_n$ if $d^-_\ell$ is not a maximum;
\item  $a^{t^+_{\ell,d,n}}_\ell=d^+_\ell$ if $d^+_\ell$ is minimum;
\item $[d^+_\ell, a^{t^+_{\ell,d,n}}_\ell)\subseteq U^\ell_n$ if $d^+_\ell$ is not a minimum.
\end{itemize} 

\medskip

For each $d\in D$, and $n\in\omega$ let $s^n_d,$ be equal to the immediate predecessor of $d$ in $C$ if it exists, and otherwise let $s^n_t=\max\{\, t^-_{\ell,d,n}\colon \ell<2m\,\}$. 

\medskip

Let $t^n_d$ to be equal to the immediate successor of $d$ in $C$ if it exists, and otherwise set $t^n_d= \min\{\, t^+_{\ell,d,n}\colon \ell<2m\,\}$. Set $U_n=\bigcup_{d\in D} (s^n_d,t^n_d)$ and let $G=\bigcap_{n\in\omega}$. It follows from the construction that $G\cap (C\cup D)=D$ as required. This concludes the proof of the Lemma. \endproof    
 \section{Open Problems}
It might  be interesting to carry a deeper analysis on the relation between   the cardinal invariant $\hia$ and  the other classical cardinal invariants. 
In particular, the following natural questions are left open from our work.

\begin{prob}
Is $\hia=\min\{\d, non(\mathcal M)\}$?
\end{prob} 
 A natural model to test this problem  is the finite support iteration of the eventually different forcing.
 
 \medskip
 
 The following problem is also open.
 \begin{prob}
 Are the cardinals $\mathfrak s$ and $\hia$ comparable?
 \end{prob}



\begin{thebibliography}{99}
\addcontentsline{toc}{chapter}{Bibliography}


\bibitem{Be} Bekkali, M. {\it Pseudo tree algebras}. Notre Dame J. Formal Logic 42, 101--108 (2001).
\bibitem{Be1} Bekkali, M. {\it Open problems in Boolean algebras over partially ordered sets.} BLAST
2010, June 2-6, 2010, University of Colorado, Boulder. 
\bibitem{BeBo} Bekkali, Mohamed; Bonnet, Robert \emph{Rigid Boolean algebras}. Handbook of Boolean algebras, Vol. 2, 637--678, North-Holland, Amsterdam, 1989.
\bibitem{BeTo}  Bekkali, Mohamed; Todor\v{c}evi\'c, Stevo {\it Algebras that are hereditarily interval}. Algebra Universalis 73 (2015), no. 1, 87--95.
\bibitem{Blass} Blass, Andreas Combinatorial cardinal characteristics of the continuum. \emph{Handbook of set theory}. Vols. 1, 2, 3, 395--489, Springer, Dordrecht, 2010.
\bibitem {Br95} Brendle, J\"org, {\it Evasion and prediction - the Specker phenomenon and Gross spaces}, Forum Mathematicum, 1995, vol 7, issue 5, 513--542. 
\bibitem{He} Heindorf, L. {\it On subalgebras of Boolean interval algebras}. Proc. Amer. Math. Soc. 125, 2265--2274 (1997).
\bibitem{Jech}  Jech, T. \emph{Set  theory}. Springer Monographs in Mathematics. Springer-Verlag, Berlin, 2003. The third millennium edition, revised and expanded.
\bibitem{KoMo} Koppelberg, S., Monk, J.D.: \emph{Pseudo-trees and Boolean algebras}. Order 8, 359--374 (1991/92).
\bibitem{Kunen} Kunen, K. \emph{An introduction to independence proofs}, volume 102 of Studies in Logic and the Foundations of Mathematics. North-Holland, 1983.
\bibitem{Mi} Miller, Arnold W. \emph{Special subsets of the real line}. Handbook of set-theoretic topology, 201--233, North-Holland, Amsterdam, 1984. 
\bibitem{MB} Monk, J.D., Bonnet, R. (eds.): Handbook of Boolean algebras, vol.1. North-Holland, Amsterdam (1989).  
\bibitem{MT} Mostowski, A., Tarski, A., {\it Boolesche Ringe mit geordneter Basis}. Fund. Math. 32, 69--86 (1939) (German).
\bibitem{Nikiel1}  Nikiel, J. {\it Orderability properties of a zero-dimensional space which is a continuous image of an ordered compactum}. Topology Appl. 31, 269--276 (1989).
\bibitem{Nikiel2} Nikiel, J. {\it On continuous images of arcs and compact orderable spaces}. Topology Proc. 14, 163--193 (1989).
\bibitem{Nikiel}  Nikiel, J.; Purisch, S.; Treybig, L. B. \emph{Separable zero-dimensional spaces which are continuous images of ordered compacta}. Houston J. Math. 24 (1998), no. 1, 45--56.
\bibitem{Od} Odintsov, A.A. {\it Separable images of ordered compacta}. Vestnik Moskov. Univ. Ser. I Mat. Mekh. 3, 35--38 (1989).
\bibitem{Stevo} Todor\v{c}evi\'c, S. \emph{Trees and linearly ordered sets}. In Handbook of set-theoretic topology, pages 235--293. North-Holland, Amsterdam, 1984.
\bibitem{Stevo1} Todor\v{c}evi\'c, Stevo \emph{Analytic gaps}. Fund. Math. 150 (1996), no. 1, 55--66.
\end{thebibliography}
\end{document}